\documentclass[11pt]{article}
\usepackage{amsthm}
\usepackage{epsfig}
\usepackage{caption}
\usepackage{algorithm}
\usepackage{algorithmic}
\usepackage{latexsym}
\usepackage{amsmath}
\usepackage{amsfonts}
\usepackage{graphicx}
\usepackage[normalem]{ulem}
\usepackage{array}
\usepackage{enumerate}
\usepackage{multirow}
\usepackage{color}\usepackage[dvipsnames]{xcolor}
\usepackage{tikz}\usetikzlibrary{patterns}
\usepackage[
  bookmarks=false, 
  colorlinks,
  citecolor=brown!70!black,
  linkcolor=brown!80!black,
  urlcolor=blue!70!black,
]{hyperref}
\usepackage{relsize,exscale}

\newtheorem{cor}{Corollary}
\newtheorem{thm}{Theorem}

\newcommand{\comm}[1]{#1}

\title{Descent distribution on Catalan words avoiding a pattern of length at most three}
\author{Jean-Luc {Baril}, Sergey Kirgizov and Vincent Vajnovszki\\  
LE2I, Universit\'e de Bourgogne Franche-Comt\'e\\
         B.P. 47 870, 21078 DIJON-Cedex France      \\
        {\tt e-mail:\{barjl,sergey.kirgizov,vvajnov\}@u-bourgogne.fr}}

\begin{document}
\maketitle

\begin{abstract} 
Catalan words are particular growth-restricted words over the set of non-negative integers, 
and they represent still another combinatorial class counted by the
Catalan numbers.
We study the distribution of descents on the sets of Catalan words avoiding a pattern of length at most three: for each such a pattern $p$ we provide a bivariate generating function where the coefficient of $x^ny^k$ in its series expansion is the number of length $n$ Catalan words with $k$ descents and avoiding $p$.
As a byproduct,
we  enumerate the set of Catalan words avoiding $p$, and we provide the popularity of descents on this set.
Some of the obtained enumerating sequences are not yet recorded in 
the On-line Encyclopedia of Integer Sequences.
\end{abstract}
\section{Introduction and notation}

Combinatorial objects counted by the Catalan numbers are very classical in combinatorics,
with a variety of applications in, among others, Biology, Chemistry, and Physics.
A length $n$ {\it Catalan word} is a word $w_1w_2\ldots w_n$ over the set of non-negative integers with $w_1=0$, and 
$$0\leq w_i\leq w_{i-1}+1,$$
for $i=2,3,\ldots n$.
We denote by $\mathcal{C}_n$ the set of length $n$ Catalan words, and $\mathcal{C}=\cup_{n\geq 0}\mathcal{C}_n$. 
For example, $\mathcal{C}_2=\{00,01\}$ and
$\mathcal{C}_3=\{000,001,010,011,012\}$.
It is well known that the cardinality of $\mathcal{C}_n$ is given by
 the $n$th Catalan number $\frac{1}{n+1} {{2n}\choose{n}}$, see for instance
\cite[exercise 6.19.$u$, p. 222]{Stanley99},
which is the general term   of the sequence \href{https://oeis.org/A000108}{A000108} in the On-line Encyclopedia of Integer Sequences (OEIS) \cite{Sloa}.
See also \cite{Mansour-Vajnovszki} where Catalan words are considered in the context
of the exhaustive generation of Gray codes for growth-restricted words.

A {\it pattern} $p$ is a word satisfying the property
that if $x$ appears in $p$, then all integers in the interval $[0,x-1]$ also appear in $p$. We say that a word $w_1w_2\ldots w_n$ contains the pattern $p=p_1\ldots p_k$ if there is a subsequence $w_{i_1}w_{i_2}\ldots w_{i_k}$ of $w$, $i_1<i_2< \cdots < i_k$, which is order-isomorphic to $p$. For example, the Catalan word $01012312301$ contains seven occurrences of the pattern $110$ and four occurrences of the pattern $210$. A word {\it avoids} the pattern $p$ whenever it does not contain any occurrence of $p$. We denote by $\mathcal{C}_n(p)$ the set of 
length $n$ Catalan words avoiding the pattern $p$, and $\mathcal{C}(p)=\cup_{n\geq 0}\mathcal{C}_n(p)$. For instance, $\mathcal{C}_4(012)=\{
0 0 0 0,
0 0 0 1,
0 0 1 0,
0 0 1 1,
0 1 0 0,
0 1 0 1,
0 1 1 0,
0 1 1 1\}$, and  $\mathcal{C}_4(101)=\{
0 0 0 0,
 0 0 0 1,$
 $
 0 0 1 0,
 0 0 1 1,
 0 0 1 2,
 0 1 0 0,
 0 1 1 0,
 0 1 1 1,
 0 1 1 2,
 0 1 2 0,
0 1 2 1,
 0 1 2 2,
 0 1 2 3
\}$.
For a set of words, the {\it popularity} of a pattern $p$ is the overall number of 
occurrences of $p$ within all words of the set, see 
\cite{Bona2012} where this notion was introduced, and \cite{AHP2015,Homberger,Rudo,Bkv2017} for some related results.

A {\it descent} in a word $w=w_1w_2\ldots w_n$ is an occurrence $w_iw_{i+1}$ such that $w_i>w_{i+1}$. Alternatively, a descent is an occurrence of the {\it consecutive} pattern $10$
({\it i.e.,} the entries corresponding to an occurrence of $10$ are required to be adjacent). We denote by $d(w)$ the number of descents of $w$, thus the popularity of descents on a set $S$ of words is 
$\sum_{w\in S}d(w)$.
The distribution of the number of descents has been widely studied on several classes of combinatorial objects such as permutations and words, 
since descents have some particular interpretations in fields as Coxeter groups or theory of lattice paths \cite{Ber,Ges}.

The main goal of this paper is to study the descent distribution on Catalan words
(see Table \ref{tab1} for some numerical values).
More specifically, for each pattern $p$ of length at most three, 
we give the distribution of descents on the sets 
$\mathcal{C}_n(p)$ of length $n$ Catalan words avoiding $p$.
We denote by $C_p(x,y)=\sum_{n,k\geq 0} c_{n,k}x^ny^k$ the bivariate generating function for the cardinality of words in $\mathcal{C}_n(p)$ with $k$ descents.
Plugging $y=1$
\begin{itemize}
\item[$-$] into $C_p(x,y)$, we deduce the generating function $C_p(x)$ for the set $\mathcal{C}_n(p)$, and
\item[$-$] into $\frac{\partial C_p(x,y)}{\partial y }$, we deduce the generating function for the popularity of descents in $\mathcal{C}_n(p)$.
\end{itemize}

From the definition at the beginning of this section it follows that 
a Catalan word is either the empty word, or it can uniquely be written as 
$0 (w'+1)w''$, where 
$w'$ and $w''$ are in turn Catalan words, and $w'+1$ is obtained from $w'$ by adding one 
to each of its entries.
We call this recursive decomposition {\it first return decomposition} of a Catalan word,
and it will be crucial in our further study.
It follows that $C(x)$, the generating function for the cardinality of $\mathcal C_n$,
satisfies: 
$$
C(x)=1+x\cdot C^2(x),
$$
which corresponds precisely to the sequence of Catalan numbers.

We conclude this section by explaining how Catalan words are naturally related to two
classical combinatorial classes counted by the Catalan numbers.

\subsubsection*{Catalan words vs. Dyck words}

A {\it Dyck word} is a word over $\{u,d\}$ with the same number of $u$'s and $d$'s, and
with the property that all of its prefixes contain no more $d$'s than $u$'s.
Alternatively, a Dyck word can be represented as a lattice path starting at $(0,0)$, 
ending at $(2n,0)$, and never going below the $x$-axis, consisting of 
up steps $u=(1,1)$ and down steps $d=(1,-1)$. 
There is a direct bijection
$\delta\mapsto w$ between the set of Dyck words of semilength $n$
 and $\mathcal{C}_n$: 
the Catalan word $w$ is the sequence of the lowest ordinate of the up
steps $u$ in the Dyck word $\delta$, in lattice path representation. 
For instance, the image through this bijection of the Dyck word 
$u d uu d uu dd uu ddd$ of semilength $7$ is
$0   01   12    12\in \mathcal{C}_7$.
Note that the above bijection gives a one-to-one correspondence between 
occurrences of 
the consecutive pattern $ddu$
in Dyck words and descents in Catalan words.

\subsubsection*{Catalan words vs. binary trees}

In \cite{Makinen} the author introduced an integer sequence representation for binary trees,
called {\it left-distance sequence}.
For a binary tree $T$, let consider the following labeling of its nodes: the root is labeled by $0$, a left child by the label of its parent,
and a right child by the label of its parent, plus one.
The left-distance sequence of $T$ is obtained by covering $T$ in inorder 
({\it i.e.}, visit recursively the left subtree, the root and then the right subtree of $T$)
and collecting the labels of the nodes.
In \cite{Makinen} it is showed that, for a given length, the set of 
left-distance sequences is precisely that of same length Catalan words.
Moreover, the induced bijection between Catalan words and binary trees
gives a one-to-one correspondence between descents in Catalan words and particular nodes
(left-child nodes having a right child) in binary trees.

\medskip

The remainder of the paper is organized as follows. 
In Section 2, we study the distribution of descents on the set $\mathcal{C}$ of Catalan words. As a byproduct, we deduce the popularity of descents in $\mathcal{C}$. We consider
also similar results for the obvious cases of Catalan words  avoiding a pattern of length two. In Section 3, we study the distribution and the popularity of descents on Catalan words avoiding each pattern of length three.

\section{The sets $\mathcal{C}$ and $\mathcal{C}(p)$ for $p\in\{00,01,10\}$}

Here we consider both unrestricted Catalan words and those avoiding a length two pattern.
We denote by $C(x,y)$ the bivariate generating function 
where the coefficient of $x^ny^k$ of its series expansion is the number of length $n$ 
Catalan words with $k$ descents. When we restrict to Catalan words
avoiding the pattern $p$, the corresponding generating function is denoted by $C_p(x,y)$.

\begin{thm}\label{th} 
We have 
$$C(x,y)=\frac {1-2x+2xy-\sqrt {1-4x+4x^2-4x^2y}}{2xy}.$$
\end{thm}
\comm{\proof
Let $w=0(w'+1)w''$ be the first return decomposition of 
a non-empty Catalan word $w$ with $w',w''\in\mathcal{C}$.  If $w'$ (resp. $w''$) is empty then the number $d(w)$ of descents in $w$ is the same as that of $w''$ (resp. $w'$); otherwise, we have $d(w)=d(w')+d(w'')+1$  
since there is a descent between $w'+1$ and $w''$. So, we obtain the functional equation $C(x,y) = 1+xC(x,y) + x(C(x,y)-1)+xy(C(x,y)-1)^2$ which gives the 
desired result.
\endproof}
As expected, $C(x)=C(x,1)=\frac {1-\sqrt {1-4x}}{2x}$ is the 
generating function for the Catalan numbers, and
$\frac{\partial C(x,y)}{\partial y }|_{y=1}$
is the generating function for the descent  popularity on $\mathcal{C}$, and we have the next corollary.

\begin{cor}\label{cor} The popularity of descents on the set
$\mathcal{C}_n$ is ${2n-2}\choose{n-3}$, and its generating function is
$\frac {1-4x+2x^2-(1-2x)\sqrt {1-4x}}{2x
\sqrt {1-4x}}$
(sequence \href{https://oeis.org/A002694}{A002694} in \cite{Sloa}).
\end{cor}

\begin{table}[h]
\begin{center}
\scalebox{1}{$\begin{array}{c|cccccccccc}
k\backslash n & 1 & 2 & 3 & 4 & 5 & 6 & 7 & 8&9&10\\
\hline
0 & {\bf 1}& 2& 4 & 8 & 16 & 32 & 64 & 128&256&512\\
1 &  &  & {\bf 1}& 6 & 24 & 80 & 240 & 672 &1792&4608\\
2 &  &  &  &  & {\bf 2} & 20 & 120 & 560 &2240&8064\\
3 &  &  &  &  &   &   & {\bf 5} & 70 &560&3360\\
4 &  &  &  &  &   &   &   &   &  {\bf 14} &252\\
\hline
\sum & 1 & 2 & 5 & 14 & 42 & 132 & 429 & 1430&4862&16796 \\
\end{array}$
}
\end{center}\caption{Number $c_{n,k}$ of length $n$ Catalan words  with $k$ descents for $1\leq n\leq 10$ and $0\leq k\leq 4$.}
\label{tab1}\end{table}

Catalan words of odd lengths encompass a smaller size Catalan structure. This result is stated in the next corollary, see the bold entries in Table~\ref{tab1}.

\begin{cor}\label{corr} Catalan words of length $2n+1$ with $n$ descents are enumerated by the $n${\rm th} Catalan number $\frac{1}{n+1}{2n\choose n}$.
\end{cor}
\comm{\proof
Clearly, the maximal number of descents in  a length $n$ Catalan word is $\lfloor\frac{n-1}{2}\rfloor$. Let $w$ be a Catalan word of length $2n+1$ with $n$ descents. We necessarily have $w=0(w'+1)w''$ with $w',w''\neq \epsilon$, $d(w')=\lfloor\frac{|w'|-1}{2}\rfloor$, $d(w'')=\lfloor\frac{|w''|-1}{2}\rfloor$ and $d(w)=d(w')+d(w'')+1$.
Since the length of $w$ is odd, $|w'|$ and $|w''|$ have the same parity. 
If $|w'|$ and  $|w''|$ are both even, then 
$d(w)=\frac{|w'|-2}{2}+\frac{|w''|-2}{2}+1=
\frac{|w'|+|w''|-2}{2}
<\lfloor\frac{(|w'|+|w''|+1)-1}{2}\rfloor=
\lfloor\frac{n-1}{2}\rfloor$ which gives a contradiction. 
So, $|w'|$ and $|w''|$ are both odd, and  we have 
$d(w)=\frac{|w'|-1}{2}+\frac{|w''|-1}{2}+1=\lfloor\frac{(|w'|+|w''|+1)-1}{2}\rfloor=
\lfloor\frac{n-1}{2}\rfloor$.
Thus the generating function $A(x)$ where the coefficient of $x^n$ is the number  of Catalan words of length $2n+1$ with $n$ descents satisfies $A(x)=1+xA(x)^2$ which is the generating function for the Catalan numbers.
\endproof}

There are three patterns of length two, namely $00$, $01$ and $10$, and 
Catalan words avoiding such a pattern do not have descents, thus
the corresponding bivariate generating functions collapse into one variable ones.

\begin{thm}\label{th00} For $p\in\{00,01\}$, we have  $C_{p}(x,y)=\frac{1}{1-x}$.
\end{thm}
\comm{\proof
If $p=00$ (resp. $p=01$) then $012\ldots n-1$ (resp. $0\ldots 0$) is the unique non-empty Catalan word of length $n$ avoiding $p$, and the 
statement follows.
\endproof}

\begin{thm}\label{th10} We have  $C_{10}(x,y)=\frac{1-x}{1-2x}$, 
which is the generating function for the sequence $2^{n-1}$ 
(sequence \href{https://oeis.org/A011782}{A011782} in \cite{Sloa}).
\end{thm}
\comm{\proof A non-empty Catalan word avoiding the pattern $10$ is of the form $0^k(w'+1)$ for $k\geq 1$, and with $w'\in \mathcal{C}(10)$. So, we have the functional equation 
$C_{10}(x)=1+\frac{x}{1-x}C_{10}(x)$, which gives $C_{10}(x)=\frac{1-x}{1-2x}$.
\endproof}

\section{The sets $\mathcal{C}(p)$ for a length three pattern $p$}

Here we turn our attention to patterns of length three.
There are thirteen such patterns, and we give 
the distribution and the popularity of descents on Catalan words avoiding 
each of them. Some of the obtained results are summarized in Tables \ref{Tab1} and \ref{Tab2}.

\begin{thm}\label{th012} For $p\in\{012, 001\}$,  we have
$$C_p(x,y)=\frac {1-x+x^2-x^2y}{1-2x+x^2-x^2y}.$$
\end{thm}
\comm{\proof 
A non-empty word $w\in\mathcal{C}(012)$  has its first return decomposition  $w=01^kw''$ where $k\geq 0$ and $w''\in \mathcal{C}(012)$. If $k=0$ or $w''=\epsilon$, then the number of descents in $w$ is the same as that of 
$w''$; otherwise, we have $d(w)=d(w'')+1$ (there is 
a descent between $1^k$ and $w''$). So, we obtain the functional equation $C_{012}(x,y) = 1+xC_{012}(x,y)+\frac{x^2}{1-x}+\frac{x^2}{1-x}y(C_{012}(x,y)-1)$ which gives the desired result.\\
\noindent 
A non-empty word $w\in\mathcal{C}(001)$ has the form $w=0(w'+1)0^k$ where $w'\in \mathcal{C}(001)$ and $k\geq 0$. If $k=0$ or $w'=\epsilon$, then the number of descents in $w$ 
is the same as that of $w'$; otherwise, we have $d(w)=d(w')+1$. So, we obtain the functional equation $C_{001}(x,y) = 1+x(C_{001}(x,y)-1)+\frac{x}{1-x}+\frac{x^2}{1-x}y(C_{001}(x,y)-1)$ which gives the desired result.
\endproof}

Considering the previous theorem and the 
coefficient of $x^n$ in $C_{p}(x,1)=\frac{1-x}{1-2x}$ and in
$\frac{\partial C_p(x,y)}{\partial y }|_{y=1}=\frac {x^3}{(1-2x)^2}$, we obtain the next 
corollary.

\begin{cor}\label{cor012} For $p\in\{012, 001\}$, we have $|\mathcal{C}_n(p)|=2^{n-1}$, and the popularity of descents on the set
$\mathcal{C}_n(p)$ is $(n-2)\cdot 2^{n-3}$ (sequence \href{https://oeis.org/A001787}{A001787} in \cite{Sloa}).
\end{cor}

As in the case of length two patterns, a
Catalan word avoiding $010$ does not have descents, and we have the next
theorem.

\begin{thm}\label{th010} If $p=010$, then
$C_p(x,y)=\frac{1-x}{1-2x}$ which is the generating function for the sequence 
$2^{n-1}$ (sequence \href{https://oeis.org/A011782}{A011782} in \cite{Sloa}). 
\end{thm}
\comm{\proof 
A non-empty word $w\in\mathcal{C}(010)$ 
can be written either as $w=0w'$ with $w'\in \mathcal{C}(10)$, or as $w=0(w'+1)$ with $w'\in \mathcal{C}(010)\setminus\{\epsilon\}$. So, we deduce $C_{010}(x)=1+xC_{10}(x)+x(C_{010}(x)-1)$, and the statement holds.
\endproof}

\begin{thm}\label{th021} For $p=021$, we have
 $$C_p(x,y)=\frac {1-4x+6x^2-x^2y-4x^3+3x^3y+x^4-x^4y}{(1-x)(1-2x)(1-2x+x^2-x^2y) }.$$
\end{thm}
\comm{\proof
Let $w$ be a non-empty word in $\mathcal{C}(021)$, and let $0(w'+1)w''$ its first return decomposition with $w',w''\in \mathcal{C}(021)$. Note that $w'$ belongs to $\mathcal{C}(10)$. We distinguish two cases: (1) $w'$ does not contain  $1$, and (2) otherwise.\\
In the case (1), $w'\in \mathcal{C}(01)$ ({\it i.e.}, $w'=0^k$ for some $k\geq 0$), and $w''\in \mathcal{C}(021)$. If $w'=\epsilon$ (resp. $w''=\epsilon$), then the number of descents in $w$ 
is the same as that of $w''$ (resp. $w'$); otherwise, we have $d(w)=d(w')+d(w'')+1$. 
So, this case contributes to $C_p(x,y)$ with  $xC_{01}(x,y)+x(C_{021}(x,y)-1)+xy(C_{01}(x,y)-1)(C_{021}(x,y)-1)$.

In the case (2), $w'\in \mathcal{C}(10)\setminus \mathcal{C}(01)$ and $w''\in \mathcal{C}(01)$. If $w''=\epsilon$ then $w$ and $w'$ have the same number of descents; otherwise, we have $d(w)=d(w')+d(w'')+1$. 
So, this case contributes to $C_p(x,y)$ with  
$x(C_{10}(x,y)-C_{01}(x,y))+xy(C_{10}(x,y)-C_{01}(x,y))(C_{01}(x,y)-1)$.

Taking into account these two disjoint cases, and adding the empty word, 
we deduce the functional equation $C_{021}(x,y)=1+xC_{01}(x,y)+x(C_{021}(x,y)-1)+xy(C_{01}(x,y)-1)(C_{021}(x,y)-1)+x(C_{10}(x,y)-C_{01}(x,y))+xy(C_{10}(x,y)-C_{01}(x,y))(C_{01}(x,y)-1)$, which after calculation gives  the result.\endproof}

\begin{cor}\label{cor021} For $p=021$, we have $C_p(x)=\frac {1-4x+5x^2-x^3}{(1-2x)^2(1-x)}$ 
which is the generating function for the sequence $(n-1)\cdot 2^{n-2}+1$
(sequence \href{https://oeis.org/A005183}{A005183} in \cite{Sloa}).
The popularity of descents on the set
$\mathcal{C}_n(p)$ is $(n+1)(n-2)\cdot 2^{n-5}$ with the generating function 
$\frac {x^3(1-x)}{(1-2x)^3}$ (sequence \href{https://oeis.org/A001793}{A001793} in \cite{Sloa}).
\end{cor}

\begin{thm}\label{th102} For $p\in\{102,201\}$, we have 
$$C_p(x,y)=\frac {1-3x+3x^2-2x^2y-x^3+x^3y}{ \left( 1-x
 \right)  \left( 1-3x+2x^2-2x^2y \right) }.$$
\end{thm}
\comm{\proof 
Let $w$ be a non-empty word in $\mathcal{C}(102)$, and let $0(w'+1)w''$ its first return decomposition with $w',w''\in \mathcal{C}(102)$. If $w'$ is empty, then $w=0w''$ for some $w''\in \mathcal{C}(102)$ and we have $d(w)=d(w'')$. If $w''$ is empty, then $w=0(w'+1)$ for some $w'\in \mathcal{C}(102)$ and we have $d(w)=d(w')$.  If $w'$ and $w''$ are both non-empty, then $w'\in \mathcal{C}(102)\setminus\{\epsilon\}$ and $w''\in \mathcal{C}(012)\setminus\{\epsilon\}$. We deduce the functional equation $C_{102}(x,y)=1+xC_{102}(x,y)+x(C_{102}(x,y)-1)+xy(C_{102}(x,y)-1)(C_{012}(x,y)-1)$. Finally, by Theorem \ref{th012} we obtain the desired result.

Let $w$ be a non-empty word in $\mathcal{C}(201)$, and let $0(w'+1)w''$ its first return decomposition with $w',w''\in \mathcal{C}(201)$. If $w'$ is empty, then $w=0w''$ for some $w''\in \mathcal{C}(201)$ and we have $d(w)=d(w'')$. If $w''$ is empty, then $w=0(w'+1)$ for some $w'\in \mathcal{C}(201)$ and we have $d(w)=d(w')$. If $w'$ and $w''$ are both non-empty, then $d(w)=d(w')+d(w'')+1$ and we distinguish two cases: (1) $w'$ does not contain $1$, and (2) otherwise. In the case (1), we have $w'\in \mathcal{C}(01)\setminus\{\epsilon\}$ and $w''\in \mathcal{C}(201)\setminus\{\epsilon\}$; in the case (2), $w'$ contains 1 and $w'\in \mathcal{C}(201)\setminus \mathcal{C}(01)$ and $w''\in \mathcal{C}(01)\setminus\{\epsilon\}$. 
Combining the previous cases, 
the functional equation becomes $C_{201}(x,y)=1+xC_{201}(x,y)+x(C_{201}(x,y)-1)+xy(C_{01}(x,y)-1)(C_{201}(x,y)-1)+xy(C_{201}(x,y)-C_{01}(x,y))(C_{01}(x,y)-1)$, which gives the desired result.
\endproof}

\begin{cor}\label{cor102}For $p\in\{102,201\}$, 
we have $C_p(x)=\frac {1-3x+x^2}{(1-x)(1-3x)}$ which is the generating function 
of the sequence $\frac{3^{n-1}+1}{2}$ (sequence \href{https://oeis.org/A007051}{A007051} in \cite{Sloa}).
The popularity of descents on the set
$\mathcal{C}_n(p)$ is
$(n-2)\cdot 3^{n-3}$ with the generating function $\frac{x^3}{(1-3x)^2}
$ (sequence \href{https://oeis.org/A027471}{A027471} in \cite{Sloa}).
\end{cor}

\begin{thm}\label{th120} For $p\in\{120,101\}$, we have 
$$C_p(x,y)=\frac {1-2x+x^2-x^2y}{1-3x+2x^2-x^2y}.$$
\end{thm}
\comm{\proof
Let $w$ be a non-empty word in $\mathcal{C}(120)$, and let $0(w'+1)w''$ be its first return decomposition where $w',w''\in \mathcal{C}(120)$. If $w''$ is empty, then $w=0(w'+1)$ for some $w'\in \mathcal{C}(120)$ and we have $d(w)=d(w')$; if $w'$ is empty, then $w=0w''$ for some $w''\in \mathcal{C}(120)$ and we have $d(w)=d(w')$; if $w'$ and $w''$ are not empty, then $w'\in \mathcal{C}(01)\setminus\{\epsilon\}$, $w''\in \mathcal{C}(120)\setminus\{\epsilon\}$ and $d(w)=d(w')+d(w'')+1$. We deduce the functional equation $C_{120}(x,y)=1+xC_{120}(x,y)++x(C_{120}(x,y)-1)+xy(C_{01}(x,y)-1)(C_{120}(x,y)-1)$ which gives the result.

Let $w$ be a non-empty word in $\mathcal{C}(101)$, and let $0(w'+1)w''$ be its first return decomposition where $w',w''\in \mathcal{C}(101)$. If $w'$ is empty, then $w=0w''$ for some $w''\in \mathcal{C}(101)$ and $d(w)=d(w'')$; if $w''$ is empty, then $w=0(w'+1)$ for some $w'\in \mathcal{C}(101)$ and $d(w)=d(w'')$; if $w'$ and $w''$ are not empty, then $w'\in \mathcal{C}(101)\setminus\{\epsilon\}$ and $w''\in \mathcal{C}(01)\setminus\{\epsilon\}$ and $d(w)=d(w')+d(w'')+1$. We deduce the functional equation $C_{101}(x,y)=1+xC_{101}(x,y)+x(C_{101}(x,y)-1)+ xy(C_{101}(x,y)-1)(C_{01}(x,y)-1)$ which gives the result.
\endproof}

\begin{cor}\label{cor120} For $p\in\{120,101\}$, we have $C_p(x)=\frac {1-2x}{1-3x+x^2}$ and the coefficient of $x^n$  in its  series expansion  is the $(2n-1)$th term of the Fibonacci sequence (see \href{https://oeis.org/A001519}{A001519} in \cite{Sloa}).
The popularity of descents on the set
$\mathcal{C}_n(p)$ is given by $\sum_{k=1}^{n-2} k\cdot {{n+k-2}\choose{2k}}$ which is the 
coefficient of $x^n$ in the series expansion of 
$\frac {x^3 \left( 1-x \right) }{ \left(1-3x+x^2\right)^2}$
(sequence \href{https://oeis.org/A001870}{A001870} in \cite{Sloa}).
\end{cor}

\begin{thm}\label{th011} For $p=011$, we have 
$$C_p(x,y)=\frac {1-2x+2x^2-x^3+x^3y}{ \left( 1-x \right)^3}.$$
\end{thm}
\comm{\proof Let $w$ be a non-empty word in $\mathcal{C}(011)$, and let $0(w'+1)w''$ its first return decomposition where $w',w''\in \mathcal{C}(011)$. If $w'$ (resp. $w''$) is empty, then we have $d(w)=d(w'')$ (resp. $d(w)=d(w')$);  if $w'$ and $w''$ are non-empty, then $w'\in \mathcal{C}(00)\setminus \{\epsilon\}$ and $w''\in \mathcal{C}(01)\setminus\{\epsilon\}$. We deduce the functional equation $C_{011}(x,y)=1+xC_{011}(x,y)+x(C_{00}(x,y)-1)+xy(C_{00}(x,y)-1)(C_{01}(x,y)-1)$ which gives the result.
\endproof}

\begin{cor}\label{cor011} For $p=011$, we have $C_p(x)=\frac {1-2x+2x^2}{(1-x)^3}$ and the coefficient of $x^n$  in its series expansion  is $1+{{n}\choose {2}}$ (sequence 
\href{https://oeis.org/A000124}{A000124} in \cite{Sloa}).
The popularity of descents on the set
$\mathcal{C}_n(p)$ is given by $\frac{(n-1)(n-2)}{2}$ which is the 
coefficient of $x^n$ in the series expansion of 
$\frac {x^3}{(1-x)^3}$ (sequence  \href{https://oeis.org/A000217}{A000217} in \cite{Sloa}).
\end{cor}

\begin{thm}\label{th000} For $p=000$, we have 
$$C_p(x,y)=\frac {1-x^2-x^2y}{1-x-2x^2-x^2y+x^3+x^4-x^4y}.$$
\end{thm}
\comm{\proof Let $w$ be a non-empty word in $\mathcal{C}(000)$, and let $0(w'+1)w''$ its first return decomposition where $w',w''\in \mathcal{C}(000)$. We distinguish two cases: (1) $w''$ is empty, and (2) otherwise.

In the case (1), we have $w=0(w'+1)$ for some $w'\in \mathcal{C}(000)$ and $d(w)=d(w')$. So, the generating function $A(x,y)$ for the Catalan words in this case is $A(x,y)=x C_{000}(x,y)$.

In the case (2), we set $w''=0(w'''+1)$ for some $w'''\in \mathcal{C}(000)$ and we have $w=0(w'+1)0(w'''+1)$. 
We distinguish three sub-cases: (2.a) $w'$ is empty, (2.b) $w'$ is non-empty and $w'''$ is empty, 
and (2.c) $w'$ and $w'''$ are both non-empty.

In the case (2.a), we have $w=00(w'''+1)$ with $w'''\in \mathcal{C}(000)$.  
So, the generating function for the Catalan words belonging to this 
case is $B_a(x,y)=x^2 C_{000}(x,y)$.

In the case (2.b), we have $w=0 (w'+1)0$ with $w'\in \mathcal{C}(000)\setminus\{\epsilon\}$.  So, the generating function for the corresponding Catalan words is $B_b(x,y)=x^2y (C_{000}(x,y)-1)$.

In the case (2.c), we have $w=0(w'+1)0(w'''+1)$ where $w'$ and $w'''$ are non-empty Catalan words such that $w'w'''$ is a Catalan word lying in the case (2). If $w'=0$, then $d(w'w''')=d(w''')=d(w)-1$; if $w'\neq 0$, then $d(w'w''')=d(w')+d(w''')+1=d(w)$. So, the generating function for the corresponding 
Catalan words is $B_c(x,y)=x^2yB_a(x,y)+ x^2 (B_b(x,y)+ B_c(x,y))$.

Considering $C_{000}(x,y)=1+A(x,y)+B_a(x,y)+B_b(x,y)+B_c(x,y)$, the obtained functional equations give the result.
\endproof}

\begin{cor}\label{cor000} For $p=000$, we have $C_p(x)={\frac{1-2x^2}{1-x-3x^2+x^3}}$ and the generating function for the popularity of descents in the sets $\mathcal{C}_n(p)$, $n\geq 0$, is
$$\frac {x^3(1-x)(1+2x)(1+x)}{( 1-x-3x^2+x^3)^2}.
$$
\end{cor}
Note that the sequences defined by the two generating functions in Corollary \ref{cor000} do not appear in  \cite{Sloa}.

\begin{thm}\label{th100} For $p=100$, we have 
$$C_p(x,y)=\frac {1-2x-x^2y+x^3}{1-3x+x^2-x^2y+2x^3}.$$
\end{thm}
\comm{\proof For $k\geq 1$, we define $\mathcal{A}_k\subset \mathcal{C}(100)$ as the set of Catalan words avoiding $100$ with exactly $k$ zeros, and let $A_k(x,y)$ be the generating function for $\mathcal{A}_k$.

A Catalan word $w\in\mathcal{A}_1$ is of the form $w=0(w'+1)$ with $w'\in \mathcal{C}(100)$. Since we have $d(w)=d(w')$, the generating function $A_1(x,y)$ for these words satisfies $A_1(x,y)=xC_{100}(x,y)$.

A Catalan word $w\in\mathcal{A}_k$, $k\geq 3$, is of the form $w=0^{k-2}w'$ with $w'\in \mathcal{A}_2$. Since we have $d(w)=d(w')$, the generating function $A_k(x,y)$ for these words satisfies $A_k(x,y)=x^{k-2}A_2(x,y)$.

A Catalan word $w\in\mathcal{A}_2$ has one of the three following forms:

 (1) $w=00(w'+1)$ with $w'\in \mathcal{C}(100)$; we have  $d(w)=d(w')$, and the generating function for these Catalan words is $x^2C_{100}(x,y)$.

 (2) $w=0(w'+1)0$ with $w'\in \mathcal{C}(100)\setminus\{\epsilon\}$; we have $d(w)=d(w')+1$, and  the generating function for these Catalan words is $x^2y(C_{100}(x,y)-1)$.

 (3) $w=0(w'+1)0(w''+1)$ where  $w'$ and $w''$ are non-empty and $w'w''\in \mathcal{A}_k$ for some
$k\geq 2$ ({\it i.e.}, $w'w''=0^{k-2}0(u+1)0(v+1)$ with $0(u+1)0(v+1)\in \mathcal{A}_2$). So, there are $(k-1)$ possible choices for $w'$, namely $0, 0^2, \ldots, 0^{k-2},$ and $0^{k-2}0(u+1)$. 
If $w'=0, 0^2, \ldots, 0^{k-2}$, then $d(w)=d(0(u+1)0(v+1))+1$; 
if $w'=0^{k-2}0(u+1)$ and $u\neq \epsilon$, then $d(w)=d(0(u+1)0(v+1))$; 
if $w'=0^{k-2}0(u+1)$ and $u= \epsilon$, then $d(w)=d(0(u+1)0(v+1))+1$. So, the generating function  for these words is 
$x^2yA_2(x,y)\sum_{k\geq 2} (k-2)x^{k-2} +x^2(A_2(x,y)-x^2C_{100}(x,y))\sum_{k\geq 2} x^{k-2}+x^4yC_{100}(x,y)\sum_{k\geq 2} x^{k-2}$, which is 
$\frac{x^3y}{(1-x)^2}A_2(x,y)+\frac{x^2}{1-x}A_2(x,y)+\frac{x^4y-x^2}{1-x}C_{100}(x,y)$.

Taking into account all previous cases, we obtain the following functional equations:
 \begin{enumerate}

\item[(i)] $A_1(x,y)=x C_{100}(x,y),$

\item[(ii)] $A_2(x,y)= x^2C_{100}(x,y)+ x^2y(C_{100}(x,y)-1)+\frac{x^3y}{(1-x)^2}A_2(x,y)+\frac{x^2}{1-x}A_2(x,y)+\frac{x^4y-x^2}{1-x}C_{100}(x,y),$

\item[(iii)] $A_k(x,y)=x^{k-2}A_2(x,y) \mbox{ for } k\geq 3,$

\item[(iv)] $C_{100}(x,y)=1+\sum_{k\geq 1}A_k(x,y).$
\end{enumerate}

A simple calculation gives the desired result.
\endproof}

\begin{cor}\label{cor100} For $p=100$, we have 
$C_p(x)=\frac {1-2x-x^2+x^3}{1-3x+2x^3}$, which is the 
generating function for the sequence $\lceil\frac{(1+\sqrt{3})^{n+1}}{12}\rceil$
(see \href{https://oeis.org/A057960}{A057960} in \cite{Sloa}),
and the generating function for the popularity of descents in the sets $\mathcal{C}_n(p)$, $n\geq 0$, is
$$\frac {x^3(1-x-x^2)}{(1-3x+2x^3)^2}.$$
\end{cor}

\begin{thm}\label{th110} For $p=110$, we have 
$$C_p(x)={\frac {1-3x+2x^2+x^3-x^4+x^4y}
{ \left( 1-x \right)  \left(1-3x+x^2+2x^3-x^3y\right) }}.
$$

\end{thm}
\comm{\proof
Let $w$ be a non-empty word in $\mathcal{C}(110)$, and let $0(w'+1)w''$ its first return decomposition where $w',w''\in \mathcal{C}(110)$.

Then, $w$ has one of the following forms:
\begin{itemize}
\item[$-$] $w=0(w'+1)$ where $w'\in \mathcal{C}(110)$; the generating function for these words is $xC_{110}(x,y)$.
\item[$-$] $w=0w'$ where $w'\in \mathcal{C}(110)\setminus\{\epsilon\}$; the generating function for these words is $x(C_{110}(x,y)-1)$.
\item[$-$] $w=0(w'+1)w''$ with $w'\in \mathcal{C}(00)\setminus\{\epsilon\}$ and $w''\in \mathcal{C}(10)\setminus\{\epsilon\}$; the generating function for these words is $xy(C_{00}(x,y)-1)(C_{10}(x,y)-1)$.
\item[$-$] The last form is $w=0(w'+1)w''$ where $w'\in \mathcal{C}(00)\setminus\{\epsilon\}$ and $w''\notin\mathcal{C}(10)$. So, we have  $w=012\ldots k0^{a_0}1^{a_1}\ldots (k-1)^{a_k} (w'''+k-1)$ where $k\geq 1$, $a_i\geq 1$ for $0\leq i\leq k$, and $w'''\in\mathcal{C}(110)\setminus\mathcal{C}(10)$; the generating function for these words is $y\sum_{k\geq 1}\frac{x^{2k+1}}{(1-x)^k}(C_{110}(x,y)-C_{10}(x,y))$.
\end{itemize}

Combining these different cases, we deduce the functional equation:

$$\begin{array}{ll}C_{110}(x,y)=&1+xC_{110}(x,y)+x(C_{110}(x,y)-1)+xy(C_{00}(x,y)-1)(C_{10}(x,y)-1)+\\
&y\sum_{k\geq 1} \frac{x^{2k+1}}{(1-x)^k}(C_{110}(x,y)-C_{10}(x,y)).\end{array}$$

Considering Theorems \ref{th10} and \ref{th00}, the result follows. 
\endproof}

\begin{cor}\label{cor110} For $p=110$, we have $C_p(x)=\frac{1-3x+2x^2+x^3}{(1-x)^2(1-2x-x^2)}$ and the generating function for the popularity of descents in the sets $\mathcal{C}_n(p)$, $n\geq 0$, is
$$\frac {x^3(1-x-x^2)^2}
{(1-x)^3(1-2x-x^2)^2}.$$
\end{cor}

\begin{thm}\label{th210} For $p=210$, we have 
$$C_p(x)=\frac {1-5x+8x^2-x^2y-4x^3+3x^3y-x^4y}
{(1-2x)(1-4x+4x^2-x^2y+x^3y)}.$$
\end{thm}
\comm{\proof
Let $w$ be a non-empty word in $\mathcal{C}(210)$, and let $0(w'+1)w''$ be its first return decomposition where $w',w''\in \mathcal{C}(210)$.

Then, $w$ has one of the following forms:
\begin{itemize}
\item[$-$] $w=0(w'+1)$ where $w'\in \mathcal{C}(210)$; the generating function for these words is $xC_{210}(x,y)$.
\item[$-$] $w=0w''$ where $w''\in\mathcal{C}(210)\setminus\{\epsilon\}$; the generating function for these words is $x(C_{210}(x,y)-1)$.
\item[$-$] $w=0(w'+1)w''$ where $w'\in \mathcal{C}(01)\setminus\{\epsilon\}$ and $w''\in \mathcal{C}(210)\setminus\{\epsilon\}$; the generating function for these sets is $xy(C_{01}(x,y)-1)(C_{210}(x,y)-1)$.
\item[$-$] $w=01^{a_1}2^{a_2} \ldots k^{a_k}w ''$ where $k\geq 2$, $a_i\geq 1$ for $1\leq i\leq k$, and $w''\in \mathcal{C}(10)\setminus\{\epsilon\}$; the generating function for these words is $y(C_{10}(x,y)-1)\sum_{k\geq 2} \frac{x^{k+1}}{(1-x)^k}$.
\item[$-$] $w=01^{a_1}2^{a_2} \ldots k^{a_k}0^{b_0}1^{b_1} \ldots (k-2)^{b_{k-2}}(w''+k-2)$ where $k\geq 2$, $a_i\geq 1$ for $1\leq i\leq k$, $b_i\geq 1$ for $1\leq i\leq k-2$, and $w''\in \mathcal{C}(210)\setminus\mathcal{C}(10)$; the generating function for these words is $y(C_{210}(x,y)-C_{10}(x,y))\sum_{k\geq 2} \frac{x^{k+1}}{(1-x)^k}\frac{x^{k-1}}{(1-x)^{k-1}}$.
\end{itemize}

Combining these different cases, we deduce the functional equation:

$$\begin{array}{ll}C_{210}(x,y)=&1+xC_{210}(x,y)+x(C_{210}(x,y)-1)+ xy(C_{01}(x,y)-1)(C_{210}(x,y)-1)+\\&y(C_{10}(x,y)-1)\sum_{k\geq 2} \frac{x^{k+1}}{(1-x)^k}+
y(C_{210}(x,y)-C_{10}(x,y))\sum_{k\geq 2} \frac{x^{k+1}}{(1-x)^k}\frac{x^{k-1}}{(1-x)^{k-1}}.\end{array}$$

Finally, considering Theorem \ref{th10} the desired result follows.
\endproof}

\begin{cor}\label{cor210} For $p=210$, we have 
$C_p(x)=\frac {1-5x+7x^2-x^3-x^4}{(1-2x)
(1-4x+3x^2+x^3)}$ and the generating function for the popularity of descents in the set 
$\mathcal{C}_n(p)$, $n\geq 0$, is
$$\frac{x^3(1-2x)}{(1-4x+3x^2+x^3)^2}.$$
\end{cor}

\begin{table}
\begin{center}
\begin{tabular}{|c|c|c|c|}
\hline
Pattern $p$ & Sequence $|\mathcal{C}_n(p)|$  & Generating function & OEIS \cite{Sloa}  \\
\hline
$012$, $001$, $010$         & $2^{n-1}$       &  $\frac{1-x}{1-2x}$  & \href{https://oeis.org/A011782}{A011782} \\
$021$         & $(n-1)\cdot 2^{n-2}+1$  &$\frac{1-4x+5x^2-x^3}{(1-x)(1-2x)^2}$
& \href{https://oeis.org/A005183}{A005183} \\
$102$, $201$  & $\frac{3^{n-1}+1}{2}$ & $\frac {1-3x+x^2}{(1-x)(1-3x)}$ & \href{https://oeis.org/A007051}{A007051}  \\
$120$, $101$          & $F_{2n-1}$   &   $\frac{1-2x}{1-3x+x^2}$    & \href{https://oeis.org/A001519}{A001519} \\
$011$       &      $\frac{n(n-1)}{2}+1$ &$\frac{1-2x+2x^2}{(1-x)^3}$& \href{https://oeis.org/A000124}{A000124}   \\
$000$                & &$\frac {1-2x^2}{1-x-3x^2+x^3}$
     &  \\
$100$       &   $\lceil\frac{(1+\sqrt{3})^{n+1}}{12}\rceil$    & $\frac{1-2x-x^2+x^3}{1-3x+2x^3}$& \href{https://oeis.org/A057960}{A057960} \\
$110$       &  $\frac{1}{2}\,\sum _{k=0}^{\lfloor \frac{n}{2}\rfloor }{n+1
\choose 2\,k+1}{2}^{k}-\frac{n-1}{2}$
  & $\frac{1-3x+2x^2+x^3}{(1-x)^2(1-2x-x^2)} $   &  \\
$210$  &  & $\frac{1-5x+7x^2-x^3-x^4}{(1-2x)(1-4x+3x^2+x^3)}$&
\\
 \hline
\end{tabular}
\end{center}
\caption{\label{Tab1}Catalan words avoiding a pattern of length three.}
\end{table}

\begin{table}
\begin{center}
\begin{tabular}{|c|c|c|c|}
\hline
\multirow{2}{*}{Pattern $p$} & Popularity of descents  & \multirow{2}{*}{Generating function} & 
\multirow{2}{*}{OEIS \cite{Sloa}} \\
            & on $\mathcal{C}_n(p)$   &                     &   \\
\hline
$012$, $001$         & $(n-2)\cdot 2^{n-3}$       &  $\frac{x^3}{(1-2x)^2}$  & \href{https://oeis.org/A001787}{A001787} \\
$010$         & $0$       &  $0$  &  \\
$021$         & $(n+1)(n-2)\cdot 2^{n-5}$  &$\frac {x^3(1-x) }{(1-2x)^3}$
& \href{https://oeis.org/A001793}{A001793} \\
$102$, $201$ & $(n-2)\cdot 3^{n-3}$ & $\frac {x^3}{(1-3x)^2}
$ & \href{https://oeis.org/A027471}{A027471} \\
$120$, $101$          & $\sum_{k=1}^{n-2} k\cdot {{n+k-2}\choose{2k}}$&   $\frac {x^3( 1-x) }
{(1-3x+x^2)^2}$    & \href{https://oeis.org/A001870}{A001870} \\
$011$       &      $\frac{(n-1)(n-2)}{2}$ &$\frac {x^3}{(1-x)^3}$& \href{https://oeis.org/A000217}{A000217}   \\
$000$                & &$\frac {x^3(1-x)(1+2x)(1+x)}{(1-x-3x^2+x^3)^2}
$
     &  \\
$100$       &       & $\frac{x^3(1-x-x^2)}{(1-3x+ 2x^3)^2}
$& \\
$110$       &
  & $\frac {x^3(1-x-x^2)^2}{
 (1-x)^3(1-2x-x^2)^2}$   &  \\
$210$  &  & $\frac{x^3(1-2x)}{(1-4x+3x^2+x^3)^2}$&
\\
 \hline
\end{tabular}
\end{center}
\caption{\label{Tab2}Popularity of descents on Catalan words avoiding a pattern of length three.}
\end{table}


\newpage

\section{Final remarks}

At the time of writing this paper, the enumerating sequences $(|\mathcal{C}_n(p)|)_{n\geq 0}$, for 
$p\in\{000,110,210\}$, are not recorded in \cite{Sloa}, and it will be interesting to explore potential
connections of these sequences with other known ones.

According to Theorem \ref{th012}, for any $k\geq 0$, the set of fixed length Catalan words 
with $k$ descents avoiding $p=001$ is 
equinomerous with those avoiding $q=012$, and a natural question that arises is to find a constructive bijection
between the two sets; and similarly for $(p,q)=(102,201)$, see Theorem \ref{th102}, and for 
$(p,q)=(101, 120)$, see Theorem \ref{th120}.
In the same vein, some of the enumerating sequences obtained in this paper count classical combinatorial classes
(see Tables \ref{Tab1} and \ref{Tab2})
and these results deserve bijective proofs.

Finally, our initiating study on pattern avoidance on Catalan words can naturally be extended to  patterns of length more than three, vincular patterns and/or
multiple pattern avoidance.
For example, some of the patterns we considered here hide larger length patterns
(for instance, an occurrence of $210$ in a Catalan word is a part of an occurrence
of $01210$), and some of our results can be restated in this light.


\begin{thebibliography}{10}

\bibitem{AHP2015}
M.~Albert, C.~Homberger, and J.~Pantone.
\newblock Equipopularity classes in the separable permutations.
\newblock {\em The Electronic J. of Comb.}, 22(2):P2.2, 2015.
\newblock (electronic).

\bibitem{Bkv2017}
J.-L. Baril, S.~Kirgizov, and V.~Vajnovszki.
\newblock Patterns in treeshelves.
\newblock {\em Discrete Mathematics}, 340(12):2946--2954, 2017.

\bibitem{Ber}
F.~Bergeron, N.~Bergeron, R.B. Howlett, and D.E. Taylor.
\newblock A decomposition of the descent algebra of a finite {C}oxeter group.
\newblock {\em J. Algebraic Combin.}, 1:23--44, 1992.

\bibitem{Bona2012}
M.~B{\'o}na.
\newblock Surprising symmetries in objects counted by {C}atalan numbers.
\newblock {\em The Electronic J. of Comb.}, 19(1):P62, 2012.
\newblock (electronic).

\bibitem{Deu}
E.~Deutsch.
\newblock Dyck path enumeration.
\newblock {\em Discrete Math.}, 204:167--202, 1999.

\bibitem{Ges}
I.~Gessel and G.~Viennot.
\newblock Binomial determinants, paths, and hook length formulae.
\newblock {\em Adv. Math.}, 58:300--321, 1985.

\bibitem{Homberger}
C.~Homberger.
\newblock Expected patterns in permutation classe.
\newblock {\em The Electronic J. of Comb.}, 19(3):P43, 2012.
\newblock (electronic).

\bibitem{Makinen}
E.~M\"akinen.
\newblock Left distance binary tree representations.
\newblock {\em BIT Numerical Mathematics}, 27(2):163--169, 1987.

\bibitem{Mansour-Vajnovszki}
T.~Mansour and V.~Vajnovszki.
\newblock Efficient generation of restricted growth words.
\newblock {\em Information Processing Letters}, 113:613--616, 2013.

\bibitem{Rudo}
K.~Rudolph.
\newblock Pattern popularity in $132$-avoiding permutations.
\newblock {\em The Electronic J. of Comb.}, 20(1):P8, 2013.
\newblock (electronic).

\bibitem{Sloa}
N.J.A. Sloane.
\newblock The on-line encyclopedia of integer sequences.
\newblock {A}vailable electronically at {\tt http://oeis.org}.

\bibitem{Stanley99}
R.P. Stanley.
\newblock {\em Enumerative Combinatorics}, volume~2.
\newblock Cambridge University Press, 1999.

\end{thebibliography}
\end{document}